\pdfoutput=1

\documentclass[11pt]{scrartcl}
\usepackage{amsthm}

\usepackage{scalerel}
\usepackage{tikz}
\usetikzlibrary{svg.path}

\usepackage[a4paper, total={6in, 9in}]{geometry}

\usepackage{graphicx}
\usepackage{dirtytalk}
\usepackage{comment}
\usepackage{bm,amsfonts,nicefrac,latexsym,amsmath,amsfonts,amsbsy,amscd,amsxtra,amsgen,amsopn,bbm,amssymb}
\usepackage[ruled]{algorithm2e}
\usepackage{array}
\usepackage{multirow}
\usepackage[pdftex,hypertexnames=false,colorlinks]{hyperref} 
\usepackage[colorinlistoftodos]{todonotes}
\usepackage{subfig}
\usepackage{floatrow}
\usepackage{placeins}
\newcommand{\R}{\mathbb{R}}

\newcommand{\lion}{\textcolor{black}}

\usepackage{soul}
\usepackage{tikz}
\usetikzlibrary{calc}

\definecolor{orcidlogocol}{HTML}{A6CE39}
\tikzset{
	orcidlogo/.pic={
		\fill[orcidlogocol] svg{M256,128c0,70.7-57.3,128-128,128C57.3,256,0,198.7,0,128C0,57.3,57.3,0,128,0C198.7,0,256,57.3,256,128z};
		\fill[white] svg{M86.3,186.2H70.9V79.1h15.4v48.4V186.2z}
		svg{M108.9,79.1h41.6c39.6,0,57,28.3,57,53.6c0,27.5-21.5,53.6-56.8,53.6h-41.8V79.1z M124.3,172.4h24.5c34.9,0,42.9-26.5,42.9-39.7c0-21.5-13.7-39.7-43.7-39.7h-23.7V172.4z}
		svg{M88.7,56.8c0,5.5-4.5,10.1-10.1,10.1c-5.6,0-10.1-4.6-10.1-10.1c0-5.6,4.5-10.1,10.1-10.1C84.2,46.7,88.7,51.3,88.7,56.8z};
	}
}
\providecommand{\keywords}[1]
{
	\small	
	\textbf{\textit{Keywords---}} #1
}

\newcommand\orcidicon[1]{\href{https://orcid.org/#1}{\mbox{\scalerel*{
				\begin{tikzpicture}[yscale=-1,transform shape]
				\pic{orcidlogo};
				\end{tikzpicture}
			}{|}}}}
\usepackage{hyperref}



\begin{document}

\title{Model-based feature selection for neural networks: A mixed-integer programming approach\thanks{Supported by C3.ai Digital Transformation Institute, and Digital Futures at KTH}}

\author{Shudian Zhao \orcidicon{0000-0001-6352-0968} 
	\footnote{Optimization and Systems Theory, Department of Mathematics, KTH Royal Institute of Technology, Stockholm, Sweden, \href{mailto:shudian@kth.se}{shudian@kth.se}, \href{mailto:jankr@kth.se}{jankr@kth.se}} 
 \, Calvin Tsay \orcidicon{0000-0003-2848-2809} \footnote{Department of Computing, Imperial College London, London, UK, \href{mailto:c.tsay@imperial.ac.uk}{c.tsay@imperial.ac.uk}}
 \and Jan Kronqvist \orcidicon{0000-0003-0299-5745} $^\dagger$\footnote{Corresponding author}}

\maketitle              
\begin{abstract}
In this work, we develop a novel input feature selection framework for ReLU-based deep neural networks (DNNs), which builds upon a mixed-integer optimization approach. While the method is generally applicable to various classification tasks, we focus on finding input features for image classification for clarity of presentation. 
The idea is to use a trained DNN, or an ensemble of trained DNNs, to identify the salient input features. The input feature selection is formulated as a sequence of mixed-integer linear programming (MILP) problems that find sets of sparse inputs that maximize the classification confidence of each category. These ``inverse'' problems are regularized by the number of inputs selected for each category and by distribution constraints. Numerical results on the well-known MNIST and FashionMNIST datasets show that the proposed input feature selection allows us to drastically reduce the size of the input to $\sim$15\% while maintaining a good classification accuracy. This allows us to design DNNs with significantly fewer connections, reducing computational effort and producing DNNs that are more robust towards adversarial attacks. 

\end{abstract}

\section{Introduction}
Over the years, there has been an active interest in algorithms for training sparse deep neural networks (DNNs) or sparsifying trained DNNs. By sparsifying a DNN we mean removing some connections (parameters) in the network, which can be done by setting the corresponding weights to zero. Examples of algorithms for sparsifying or training sparse networks include, dropout methods \cite{gal2016dropout,hinton2012improving,kingma2015variational,wan2013regularization}, optimal/combinatorial brain surgeon \cite{hassibi1992second,yu2022combinatorial}, optimal brain damage \cite{lecun1989optimal}, and regularization based methods \cite{liu2015sparse,manngaard2018structural,wen2016learning}. Carefully sparsifying the network, \textit{i.e.,} reducing the number of parameters wisely, has shown to reduce over-fitting and improve overall generalizability \cite{hinton2012improving,labach2019survey}. This paper focuses on feature selection for DNNs, which can also be interpreted as ``sparsifying'' the first/input layer, and we show that we can significantly reduce the number of parameters, \textit{i.e.}, non-zero weights, while keeping a good accuracy. Throughout the paper we focus on image classification, but the framework is general.

This work focuses on finding the salient input features for classification using a DNN. We hypothesize that the number of inputs to DNNs for classification can often be greatly reduced by a \say{smart} choice of input features while keeping a good accuracy (\textit{i.e.,} feature selection). We build the hypothesis on the assumption that not all inputs, or pixels, will be equally important. Reducing the number of inputs/parameters has the potential to: i) reduce over-fitting, ii) give more robust DNNs that are less sensitive for adversarial attacks (fewer degrees of freedom for the attacker), and iii) reduce computational complexity both in training and evaluating the resulting classifier (fewer weights to determine and fewer computational operations to evaluate the outputs). The first two are classical focus areas within artificial intelligence (AI), and the third is becoming more important with an increasing interest in so-called green AI \cite{schwartz2020green}. 
Most strategies for feature selection can be grouped as either \textit{filter} methods, which examine the data, \textit{e.g.,} for correlation, and \textit{wrapper} methods, which amount to a guided search over candidate models \cite{chandrashekar2014survey,li2017feature}. 
Feature selection can also be incorporated directly into model training using \textit{embedded} methods, \textit{e.g.}, regularization. This paper and the numerical results are intended as a proof of concept to demonstrate mixed-integer linear programming (MILP) as an alternative technology for extracting the importance of input features from DNNs. 
Input feature selection is an active research area, \textit{e.g.,} see the review papers \cite{ghojogh2019feature,zebari2020comprehensive}, and a detailed comparison to state-of-the-art methods is not within the scope of this paper.

Our proposed method leverages trained models that achieve desirable performance, and attempts to select a feature set that replicates the performance using mixed-integer programming. 
We build on the idea that, given a relatively well-trained DNN, we can analyze the DNN in an inverse fashion to derive information about the inputs. 
Specifically, to determine the most important inputs, or pixels in the case of image classification, for a given label, we solve an optimization problem that maximizes the classification confidence of the label with a cardinality constraint on the number of non-zero inputs to the DNN. We consider this as an ``inverse problem,'' as the goal is to determine the DNN inputs from the output. 
We additionally propose some input distribution constraints to make the input feature selection less sensitive to errors in the input-output mapping of the DNN. We only consider DNNs with the rectified linear unit (ReLU) activation function, as it enables the input feature selection problem to be formulated as a MILP problem \cite{fischetti2018deep,lomuscio2017approach}. However, the framework can be easily generalized to CNN architectures and other MIP representable activation functions, e.g., max pooling and leaky ReLU

Optimizing over trained ReLU-based DNNs has been an active research topic in recent years, and has a wide variety of applications including verification \cite{botoeva2020efficient,lomuscio2017approach,tjeng2017evaluating}, lossless compression \cite{serra2020lossless}, and surrogate model optimization \cite{grimstad2019relu,yang2022modeling}. There even exists software, such as OMLT \cite{ceccon2022omlt}, for directly incorporating ReLU DNNs into general optimization models.   Optimizing over a trained ReLU-based DNN through the MILP encoding is not a trivial task, but significant progress has been made in terms of strong formulations \cite{anderson2020strong,kronqvist2021between,tsay2021partition}, solution methods \cite{de2021scaling,perakis2022optimizing}, and  techniques for deriving strong valid inequalities \cite{anderson2020strong,botoeva2020efficient}. In combination with the remarkable performance of state-of-the-art MILP solvers, optimization over DNNs appears computationally tractable (at least for moderate size DNNs). 
This work builds upon recent optimization advancements, as reliably optimizing over ReLU-based DNNs is a key component in the proposed method.  

The paper is structured as follows. Section~\ref{sec:formulation} first describes the MILP problem to determine which inputs maximize the classification confidence. Some enhancements for the input selection are presented, and the complete input selection algorithm is presented in Section~\ref{sec:formulation_numcontrol}. Numerical results are presented in Section~\ref{sec:results}, where we show that we can obtain a good accuracy \lion{when downsizing the input to ~15\%} by the proposed \lion{algorithm}, and that the resulting DNNs are more robust towards adversarial attacks in the $\ell_\infty$ sense. Section~\ref{sec:summary} provides some conclusions.

\section{Input feature selection algorithm} \label{sec:formulation}

Our feature selection strategy is based on the idea of determining a small optimal subset of inputs, or pixels for the case of image classification, that are allowed to take non-zero values to maximize the classification confidence for each label using a pre-trained DNN. 
By combining the optimal subsets for each label, we can determine a set of salient input features. These input features can be considered as the most important for the given DNN, but we note that the DNN might not offer a perfect input-output mapping. To mitigate the impact of model errors, we propose a technique of using input distribution constraints to ensure that the selected input features are to some extent distributed over the input space. This framework could easily be extended to use optimization over an ensemble DNN model \cite{wang2021acceleration} for input selection, where the inputs would be selected such that the ensemble classification confidence is maximized. While using DNN ensembles can further mitigate the effect of errors in individual DNNs, our initial tests did not indicate clear advantages of using DNN ensembles for this purpose. 

Here we focus on fully connected DNNs that classify grayscale images into 10 categories. While this setting is limited, the proposed method is applicable to classification of RGB images, and other classification problems in general.  The input features are scaled between 0 and 1, with 255 (white) corresponding to 1 and 0 remaining black. 
We start by briefly reviewing MILP encoding of DNNs in the next subsection, and continue with more details on the input feature selection in the following subsections.  

\subsection{Encoding DNNs as MILPs}

In a fully connected ReLU-based neural network, the $l$-th layer with input  $x^l$ and output $x^{l+1}$ is described as
$$x^{l+1} = \max \{0, W^lx^l + b^l\},$$
where $W^l \in \R^{n_{l+1}\times n_{l}}$ is the weight matrix and $b^l \in \R^{n_{l+1}}$ is the bias vector.

The input-output mapping of the ReLU activation function is given by a piece-wise linear function, and is mixed-integer representable \cite{vielma2015mixed}.
There are different formulations for encoding the ReLU activation function using MILP, where the big-M formulation \cite{fischetti2018deep,lomuscio2017approach} was the first presented MILP encoding and remains a common approach. For the $i$-th ReLU node at a fully connected layer with input $x^l$, the big-M formulation for input-output relation is given by
\begin{equation}
\begin{aligned}
\label{eq:big-M}
    (w_i^l)^\top x^l + b^l_i \leq x^{l+1}_i, \\
   (w^l_i)^\top x^l + b^l_i - (1-\sigma)LB^{l+1}_i \geq x^{l+1}_i,\\
    x^{l+1}_i \leq \sigma UB^{l+1}_i,\\
    \sigma \in \{0,1\},~x^{l+1}_i\geq 0,
\end{aligned}
\end{equation}
where \lion{$w^l_i$ is the $i$-th row vector of $W^l$, $b^l_i$ is the $i$-th entry of $b^l$},  $LB^{l+1}_i$ and $UB^{l+1}_i$ are upper and lower bounds on the pre-activation function over \lion{$x^{l+1}_i$}, such that $LB^{l+1}_i \leq (w_i^l)^\top x^l + b^l_i \leq UB^{l+1}_i$. 

The big-M formulation is elegant in its simplicity, but it is known to have a weak continuous relaxation which may require the exploration of a huge number of branch-and-bound nodes in order to solve the problem \cite{anderson2020strong}. Anderson et al. \cite{anderson2020strong} presented a so-called extended convex hull formulation, which gives the strongest valid convex relaxation of each individual node, and a non-extended convex hull formulation. Even though the convex hull is the strongest formulation for each individual node, it does not in general give the convex hull of the full input-output mapping of the DNN. Furthermore, the convex hull formulation results in a large problem formulation that can be computationally difficult to work with. The class of partition-based, or $P$-split, formulations, was proposed as an alternative formulation with a stronger continuous relaxation than big-M and computationally cheaper than the convex hull~\cite{kronqvist2021between,tsay2021partition}. Computational results in \cite{kronqvist2022p,tsay2021partition} show that the partition-based formulation often gives significant speed-ups compared to the big-M or convex hull formulations. Here, we do not focus on the computational efficiency of optimizing over ReLU-based DNNs, and, for the sake of clarity, we use the simpler big-M formulation \eqref{eq:big-M}. In fact, for the problems considered in this work, the computational time to solve the MILP problems did not represent a limiting factor (big-M tends to actually perform relatively well for simple optimization problems). But, alternative/stronger formulations could directly be used within our framework. 
\subsection{The optimal sparse input features (OSIF) problem}

With the MILP encoding of the DNN, we can rigorously analyze the DNN and find extreme points in the input-output mapping. Recently Kronqvist et al. \cite{kronqvist2022p} illustrated a similar optimal sparse input features (OSIF) problem. This problem aims to maximize the probability of at most $\aleph$ non-zero input features being classified with a certain label $i$ for a given trained DNN. The problem is formed by encoding the ReLU activation function for each hidden node by MILP. Instead of the softmax function, the objective function is  $x^L_i$, where $x^L$ is the output vector, thus maximizing the classification confidence of label $i$.

\begin{figure}[!htbp]
	\centering
	\subfloat[$i=0$, $\aleph=10$]{\includegraphics[width = 1.8in,trim={1cm 0.8cm 1cm 1.5cm},clip]{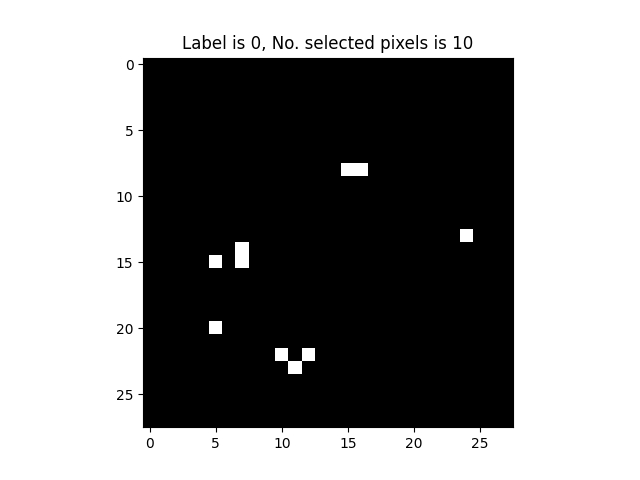}} 
	\subfloat[$i=1$, $\aleph=10$]{\includegraphics[width = 1.8in,trim={1cm 0.8cm 1cm 1.5cm},clip]{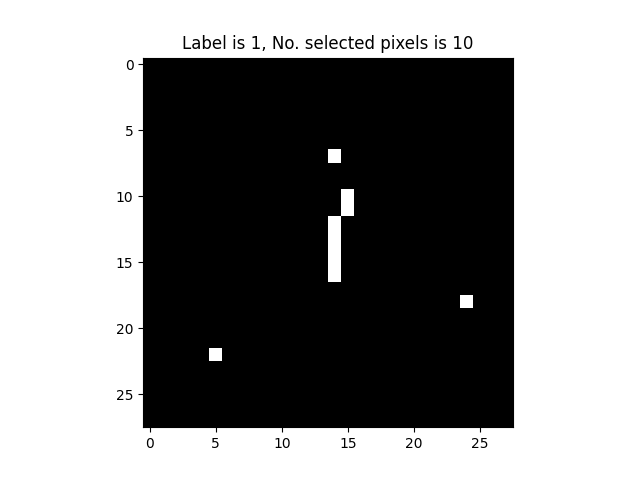}}
 
 \subfloat[$i=0$, $\aleph=20$]{\includegraphics[width = 1.8in,trim={1cm 0.8cm 1cm 1.5cm},clip]{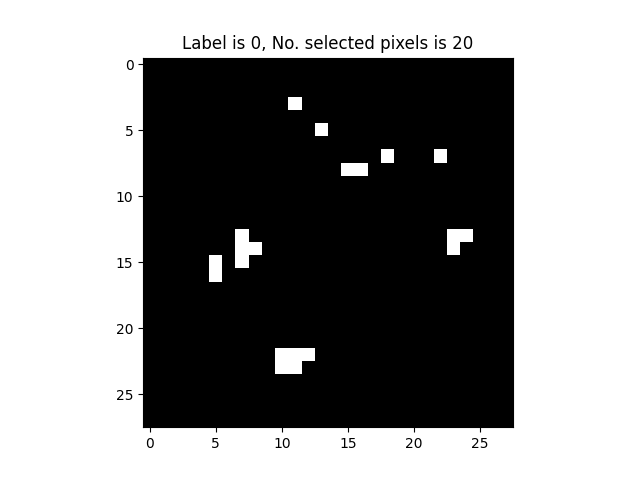}} 
	\subfloat[$i=1$, $\aleph=20$]{\includegraphics[width = 1.8in,trim={1cm 0.8cm 1cm 1.5cm},clip]{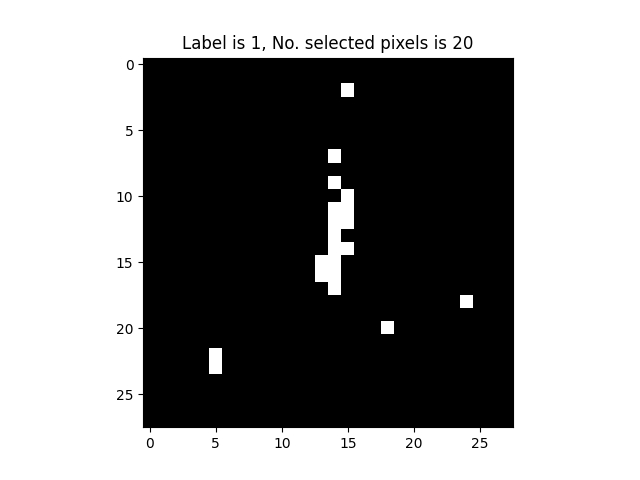}}
	\caption{The results of OSIF for class 0 and 1 on MNIST ($\aleph_0\in \{10,20\}$). Note that the selected pixels are white. }
	\label{fig:mnist_osif}. 
\end{figure}

Using the big-M formulation, the OSIF problem can be stated as
\begin{subequations}\label{eq:osif-basic}
    \begin{align}
        \max&~ x^{L}_i \\
    \textrm{s.t.}~& W^l x^l + b^l  \leq x^{l+1},~\forall l \in [L-1], \label{eq:relu_1}\\
     &W^l x^l + b^l - \textrm{diag}(LB^{l+1})(\textbf{1}-\sigma^{l+1}) \geq x^{l+1},~\forall~l \in [L-1],\label{eq:relu_2}\\
     & x^l\leq \textrm{diag}(UB^l)\sigma^l ,~\sigma^l \in \{0,1\}^{n_l}, \forall 
     ~l \in \{2,\dots ,L\},\label{eq:relu_3}\\
     & x^{L} = W^{L-1}x^{L-1}+b^{L-1}, x^{L} \in \R^{10},\\
      & x^l \in \mathbb{R}^{n_l}_+, \forall~l \in [L-1],\label{eq:relu_var_size} \\
     & y \geq x^1, y \in \{ 0, 1\}^{n_1}, \label{eq:card_1}\\
     & \mathbf{1}^\top y \leq \aleph_0, \label{eq:card_2}
    \end{align}
\end{subequations}
where $n_1 $ is the size of the input data, $x^{L} \in \R^{10}$ is the output, $L$ is the number of layers, $LB^l$ and $UB^l$ are the bounds on $x^l$, $\mathbf{1}$ denotes the all-ones vector, \lion{$\textrm{diag}(\cdot)$ denote the matrix with $\cdot$ on the diagonal and 0 on other entries}. Eq.~\eqref{eq:card_1} and \eqref{eq:card_2} describe the cardinality constraint $\| x^1\|_0 \leq \aleph_0$, which limits the number of selected inputs. 
Fig~\ref{fig:mnist_osif} shows some example results of solving problem~\eqref{eq:osif-basic} with $\aleph \in \{10,20\}$ and class $i \in \{0,1\}$. Given a larger cardinality number $\aleph$, the latter is more visually recognizable from the selected pixels.

\subsection{Input distribution constraints}
To extract information across the whole input image, we propose to add constraints to force the selected pixels to be distributed evenly across some pre-defined partitioning of the input space. Forcing the selected pixels to be more spread-out may also mitigate the effect of inaccuracy of the DNN used in the OSIF problem, \textit{e.g.,} by preventing a small area to be given overly high priority.    
There are various ways to partition the input variables. In this paper, we focus on image classification problems with square images as input. 
Furthermore, we assume that the images are roughly centered. 
Therefore, we denote each input variable as matrix $X \in \R^{n\times n}$ and define the partition as $k^2$ submatrices of equal size, \textit{i.e.,} $X^{ij} \in \R^{\frac{n}{k}\times \frac{n}{k}}$ for $i,j \in [k]$. For instance, given $n$ is even and $k=2$, a natural partition of the matrix is 
\begin{equation*}
    X=\begin{pmatrix}
       X^{11} &X^{12}\\
       X^{21} &X^{22}
    \end{pmatrix}.
\end{equation*}
In this way, we denote $x^1:=\textrm{vec}(X)$ the input data and $n_1:=n^2$ the size of the input data, then $I_{ij}$ is the index set for entries mapped from $X_{ij}$
\begin{equation*}
        I_{ij} = \{(i_1-1)n+i_2 \mid i_1 \in \{(i-1)\frac{n}{k}+1,\dots, i\frac{n}{k}-1 \}, i_2 \in \{(j-1)\frac{n}{k}+1,\dots, j\frac{n}{k}-1 \} \}.
\end{equation*}
We denote the collection of index sets for the partition as $ \mathcal{I}:=\{I_{i,j}\}_{\forall i,j \in [k]}$.
To limit the number of pixels selected from each box for each category we add the following constraints 
\begin{equation}\label{eq:spread}
\begin{aligned}
\lfloor \frac{\aleph_0}{k^2} \rfloor  \leq \sum_{i\in I_t} y_i \leq \lceil \frac{\aleph_0}{k^2} \rceil,\forall I_t \in \mathcal{I},
\end{aligned}
\end{equation}
The constraint \eqref{eq:spread} forces the pixels to spread evenly between all partitions, while allowing some to contain one more selected pixel for each category. 

\begin{figure}[!htbp]
	\centering
	\subfloat[$k=1$]{\includegraphics[width = 1.8in,trim={1cm 0.8cm 1cm 2.5cm}]{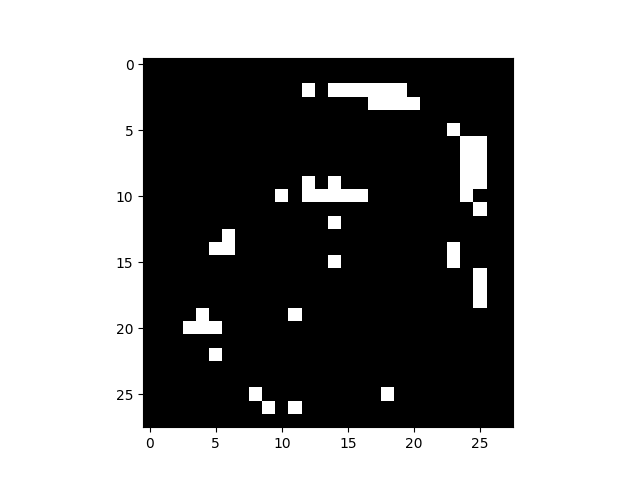}} 
	\subfloat[$k=2$]{\includegraphics[width = 1.8in,trim={1cm 0.8cm 1cm 2.5cm}]{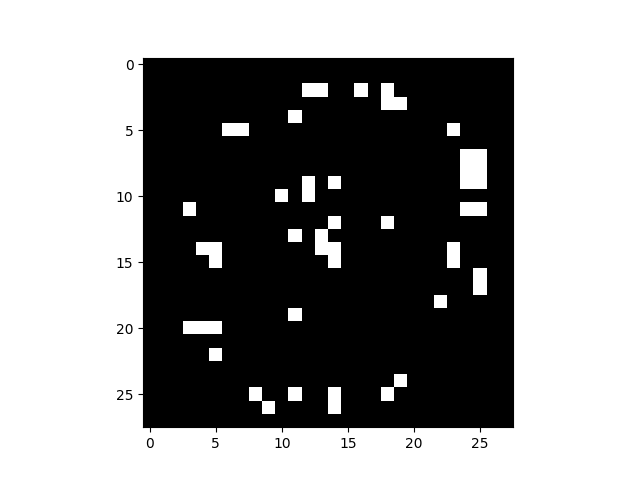}}
		
	\caption{Optimal input features for MNIST ($\aleph=50$)}
	\label{fig:mnist_50_distriConstr}
\end{figure}

To illustrate on the impact on the distribution constraints, Fig.~\ref{fig:mnist_50_distriConstr} compares the selected pixels for MNIST with $k \in \{1,2\}$. Compared to the result without distribution constraints (equivalent to $k=1$), pixels selected with $k=2$ are more scattered over the whole image and are more likely to identify generalizable distinguishing features of the full input data, assuming the dataset has been pre-processed for unused areas of the images matrices. 

\subsection{Controlling the number of selected features} \label{sec:formulation_numcontrol}
Repeatedly solving the OSIF problem \eqref{eq:osif-basic} for each label, \textit{i.e.,} $i \in \{0,\dots,9\}$, and taking the union of all selected pixels does not give us full control of the total number of selected pixels. Specifically, some of the pixels can be selected by the OSIF problem for multiple classes, resulting in fewer combined pixels (the union of selected subsets) than an initial target.

Therefore, we present an approach to control the number of selected inputs, which we use in the proposed MILP-based feature selection algorithm. The main idea is to allow freedom over features already selected by previous classes in the current OSIF problem and adjust the cardinality constraint \eqref{eq:card_2} to 
\begin{equation}\label{eq:card_numcontrol}
    \sum_{i\in [n_1] \setminus J} y_i \leq \aleph_0,
\end{equation}
where $J$ is the index set for input features selected by previous models.
Similarly, constraints \eqref{eq:spread} are adjusted as
\begin{equation}\label{eq:spread_numcontrol}
    \begin{aligned}
       \lfloor \frac{\aleph_0}{k^2} \rfloor \leq \sum_{i\in I_t\setminus J} y_i \leq \lceil \frac{\aleph_0}{k^2} \rceil,\forall I_t \in \mathcal{I}.
    \end{aligned}
\end{equation}
Finally, we formulate the OSIF problem with input distribution and total number control as 
\begin{equation}
    \begin{aligned}
    OSIF(\mathcal{M},i, \aleph_0,\mathcal{I},J) =\textrm{argmax}~ \{x^L_i~
    \mid \textrm{\eqref{eq:relu_1}--\eqref{eq:card_1}}, 
    \eqref{eq:card_numcontrol}, 
    \eqref{eq:spread_numcontrol} \}.
    \end{aligned}
\end{equation}

Based on the described techniques, we introduce the input feature selection algorithm, which is presented as pseudo code in Alg.~\ref{alg:input_select}.  
\begin{algorithm}[!htbp]
\caption{MILP-based feature selection (MILP-based selection)}\label{alg:input_select}
\KwData{the number of features $\aleph$, a trained DNN $\mathcal{M}$, a matrix partition set $\mathcal{I}$, class set $\mathcal{C}$\;}
\KwIn{$J\leftarrow \emptyset$\;}
\KwOut{Index set $J$\;}
$\aleph_0 \leftarrow \aleph/10$\;
\For{ $i \in \mathcal{C}$ }{
    $x  \leftarrow OSIF(\mathcal{M},i,\aleph_0,\mathcal{I},J)$ ; \quad\quad \quad \quad \quad \quad \quad  \# Eq. (7) \\
    $J \leftarrow  J \cup \{ s \mid x^1_s =0, s\in [n_1]\}$\;
    }
\end{algorithm}

\section {Computational results} \label{sec:results}
In this paper, we focus on image classification problems for the MNIST \cite{deng2012mnist} and FashionMNIST~\cite{xiao2017/online} datasets. Both datasets consist of a training set of 60,000 examples and a test set of 10,000 examples. Each sample image in both datasets is a $28\times 28$ grayscale image associated with labels from 10 classes. MNIST is the dataset of handwritten single digits between $0$ and $9$. FashionMNIST is a dataset of Zalando's article images with 10 categories of clothing. There is one fundamental difference between the two data sets, besides FashionMNIST being a somewhat more challenging data set for classification. In MNIST there are significantly more pixels that do not change in any of the training images, or only change in a few images, compared to FashionMNIST. The presence of such \say{dead} pixels is an important consideration for input feature selection.

Image preprocessing and training DNNs are implemented in PyTorch \cite{NEURIPS2019_9015}, and the MILP problems are modeled and solved by Gurobi through the Python API \cite{gurobi}. We trained each DNN with 2 hidden layers of the same size. 

\subsection{Accuracy of DNNs with sparse input features}

The goal is to illustrate that Alg.~\ref{alg:input_select} can successfully identify low-dimensional salient input features. We chose to focus on small DNNs, as DNN $ 2\times 20 $ can already achieve an accuracy of ~95.7\% (resp. ~86.3\%) for MNIST (resp. FashionMNIST) and larger DNNs did not give clear improvements for the input selection. For such models, the computational cost of solving the MILPs is low\footnote{On a laptop with a 10-core CPU, Gurobi can solve instances with $\aleph=100$ and a DNN of $2\times 20$ under 15 seconds. However, previous research \cite{anderson2020strong,tsay2021partition} has shown that significant speed-ups can be obtained by using a more advanced MILP approach.}.   

Table~\ref{tab:compare_k_mnist_dnn} and Table~\ref{tab:compare_k_fashionmnist_dnn} present the accuracies of DNNs with sparse input features on MNIST and FashionMNIST. It is  possible to obtain a much higher accuracy by considering a more moderate input reduction (about 0.5\% accuracy drop with 200 --300 inputs), but this defeats the idea of finding low dimensional salient features. For grayscale input image of $28\times28$, we select at most 15\% input features and present the results with $\aleph \in \{50,100 \}$.  

\subsubsection{MILP-based feature selection}

Table~\ref{tab:compare_k_mnist_dnn} compares the accuracy of classification models with different architectures, \textit{i.e., }with $2\times 20$ vs. with $2\times 40$. We select sparse input features by Alg.~\ref{alg:input_select} with OSIF models  $2\times 10$ and $2\times 20$. 
Since the distribution constraints are supposed to force at least one pixel selected in each submatrix, we select partition number $k\in\{1,2\}$ and $k \in \{1,2,3\}$ for instances with $\aleph = 50$ and $\aleph =100$ respectively. 

\begin{table}[!htb]
    \centering
    \begin{tabular}{c|rrr|rrr}
    \multicolumn{6}{c}{DNNs of $2\times20$ hidden layers}\\
        \hline
         $\aleph$             & OSIF Model & $k$ & Acc. & OSIF Model & $k$ & Acc.   \\
         \hline
         \multirow{2}{*}{50} & $2\times 10$    & 1   & 80.6\%     &  $2 \times 20$   &  1  &    80.5\%     \\
                              & $2\times 10$    & 2   &   \textbf{85.3\%}   &  $2 \times 20$   &  2  &   \textbf{86.6\%}     \\
          \hline
          \multirow{3}{*}{100} & $2\times 10$    & 1   &  88.8\%    &  $2 \times 20$   &  1  &   89.0\%     \\
                              & $2\times 10$    & 2   &   \textbf{91.2\%}    &  $2 \times 20$  &  2  &    \textbf{90.6\%}     \\
                              & $2\times 10$   & 3   &   89.3\%   &  $2 \times 20$   &  3  &     89.3\%    \\
        \hline         
           \multicolumn{6}{c}{\phantom{a}}\\
         \multicolumn{6}{c}{DNNs of $2\times40$ hidden layers}\\
           \hline
         $\aleph$             & OSIF Model & $k$ & Acc. & OSIF Model & $k$ & Acc.   \\
         \hline
          \multirow{2}{*}{50} &  $2\times 10$    & 1   &  83.1\%    &  $2 \times 20$    &  1  &     82.7\%    \\
                              &  $2\times 10$   & 2   &  \textbf{87.6\%}    &  $2 \times 20$    &  2  &      \textbf{89.2\%}  \\
          \hline
          \multirow{3}{*}{100} &  $2\times 10$   & 1   &  91.4\%    &  $2 \times 20$    &  1  &  91.4\%     \\
                              &  $2\times 10$   & 2   &   \textbf{93.4\%}   &  $2 \times 20$    &  2  &   \textbf{92.9\%}     \\
                              &  $2\times 10$   & 3   &   92.3\%   &  $2 \times 20$    &  3  &   91.8\%      \\
        \hline   
    \end{tabular}
    \caption{Accuracy of DNNs of different architectures  with sparse input features selected by Alg.~\ref{alg:input_select} on MNIST}
    \label{tab:compare_k_mnist_dnn}
\end{table}

\begin{table}[!htb]
    \centering
   \begin{tabular}{c|rrr|rrr}
   \multicolumn{6}{c}{DNNs of $2\times20$ hidden layers}\\
           \hline
         $\aleph$             & OSIF Model & $k$ & Acc. & OSIF Model & $k$ & Acc.   \\
         \hline
          \multirow{2}{*}{50} & $2\times 10$    & 1   &  76.6\%    &  $2 \times 20$  &  1  &   77.2\%      \\
                              & $2\times 10$    & 2   &   \textbf{77.0\%}   &  $2 \times 20$  &  2  &   \textbf{77.9\%}     \\
          \hline
          \multirow{3}{*}{100} & $2\times 10$    & 1   &  81.1\%     &  $2 \times 20$  &  1  &   81.3\%     \\
                              & $2\times 10$   & 2   &   \textbf{82.3\%}   &  $2 \times 20$   &  2  &   81.8\%   \\
                              & $2\times 10$    & 3   &  82.2\%    &  $2 \times 20$   &  3  &    \textbf{82.1\%}     \\
        \hline      
        \multicolumn{6}{c}{\phantom{a}}\\
         \multicolumn{6}{c}{DNNs of $2\times40$ hidden layers}\\
        \hline
         $\aleph$             & OSIF Model & $k$ & Acc. &  OSIF Model & $k$ & Acc.   \\
         \hline
    
          \multirow{2}{*}{50} & $2\times 10$   & 1   &   78.3\%   &  $2 \times 20$   &  1  &    78.6\%     \\
                              & $2\times 10$    & 2   &  \textbf{79.4\%}    &  $2 \times 20$  &  2  &    \textbf{79.6\%}    \\
          \hline
          \multirow{3}{*}{100} & $2\times 10$   & 1   &  82.4\%    &  $2 \times 20$   &  1  &    82.6\%    \\
                              & $2\times 10$    & 2   &  83.1\%    &  $2 \times 20$   &  2  &     83.2\%    \\
                              & $2\times 10$    & 3   &  \textbf{83.7\%}    &  $2 \times 20$  &  3  &    \textbf{84.0\%}     \\
        \hline                   
    \end{tabular}
    \caption{Accuracy of DNNs of different architectures with sparse input features selected by Alg.~\ref{alg:input_select} on FashionMNIST}
    \label{tab:compare_k_fashionmnist_dnn}
\end{table}

\begin{figure}[!htbp]
	\centering
	\subfloat[MNIST, $k=1$]{\includegraphics[width = 1.5in,trim={1cm 0.5cm 1cm 1.5cm},clip]{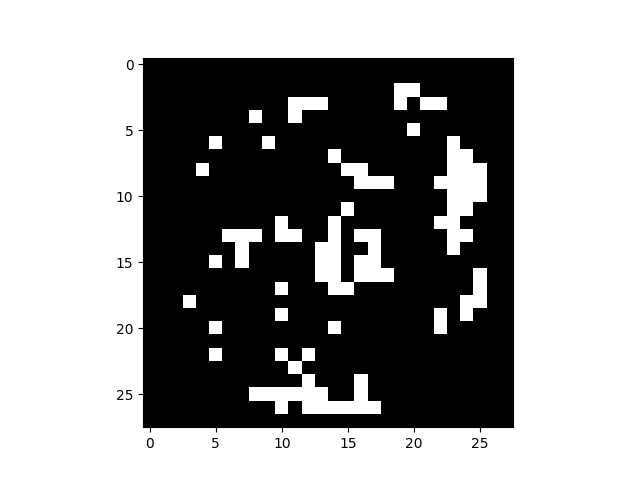} \label{fig:mnist_100_k1}} 
    \subfloat[MNIST, $k=2$]{\includegraphics[width = 1.5in,trim={1cm 0.5cm 1cm 1.5cm},clip]{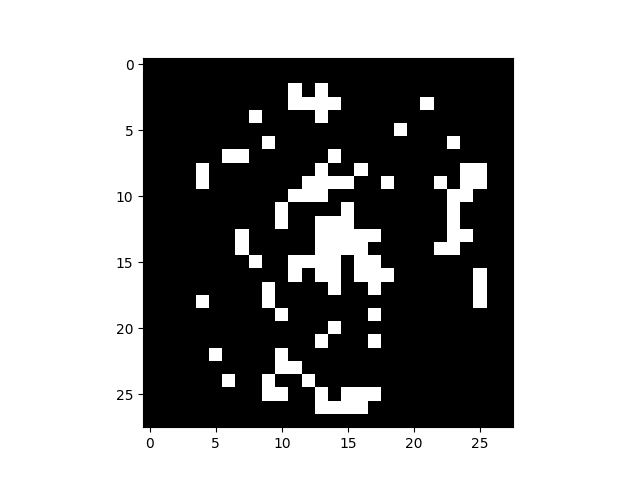}\label{fig:mnist_100_k2}}
    \subfloat[MNIST, $k=3$]{\includegraphics[width = 1.5in,trim={1cm 0.5cm 1cm 1.5cm},clip]{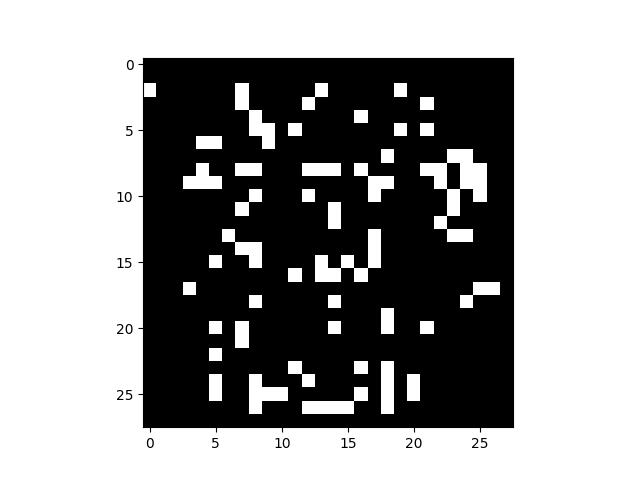}\label{fig:mnist_100_k3}}

 \subfloat[FashionMNIST, $k=1$]{\includegraphics[width = 1.5in,trim={1cm 0.5cm 1cm 1.5cm},clip]{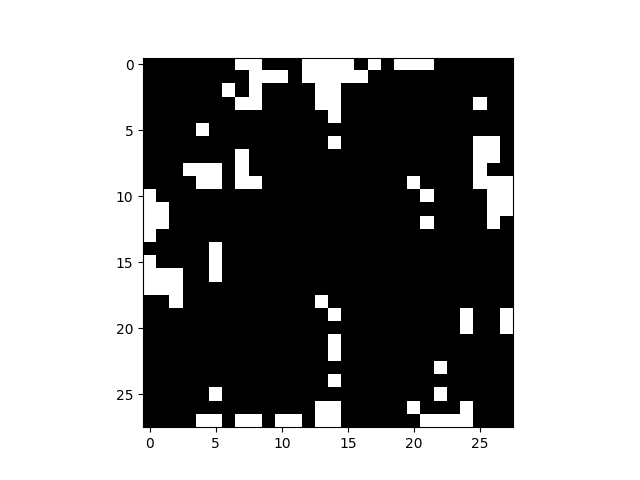}\label{fig:fashionmnist_100_k1}} 
    \subfloat[FashionMNIST, $k=2$]{\includegraphics[width = 1.5in,trim={1cm 0.5cm 1cm 1.5cm},clip]{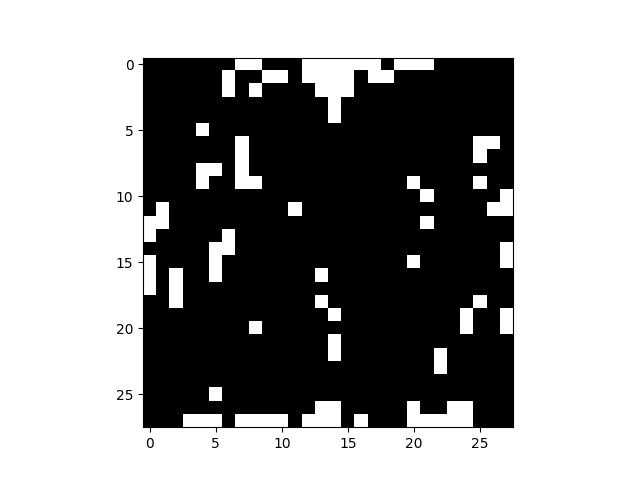}\label{fig:fashionmnist_100_k2}}
    \subfloat[FashionMNIST, $k=3$]{\includegraphics[width = 1.5in,trim={1cm 0.5cm 1cm 1.5cm},clip]{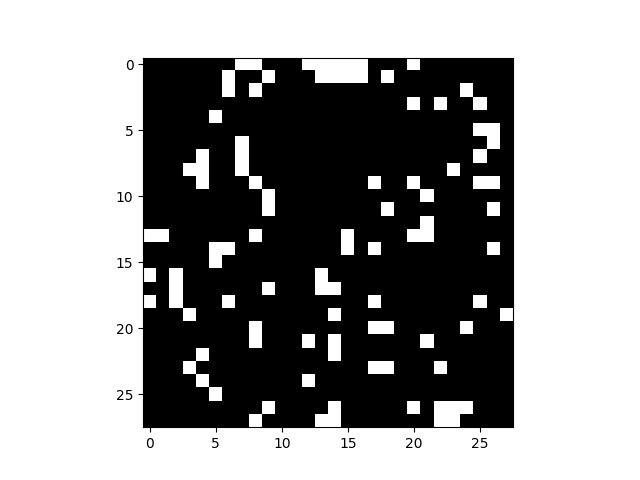}\label{fig:fashionmnist_100_k3}}
	\caption{MNIST and FashionMNIST input features selected by Alg.~\ref{alg:input_select} with $\aleph= 100$ and $k\in \{1,2,3\}$}
	\label{fig:mnist_100_distriConstr}. 
\end{figure}

First, we investigate the effect of the distribution constraints. Table~\ref{tab:compare_k_mnist_dnn} and Table~\ref{tab:compare_k_fashionmnist_dnn} both show that the accuracy increases  when adding the distribution constraints (\textit{i.e.}, from $k=1$ to $k=2$) for $\aleph \in\{50,100\}$.
However, the distribution constraints become less important as the number of selected features $\aleph$ increases; the best $k$ for MNIST and FashionMNIST varies for $\aleph =100$ (noting that the choice of $k$ also affects accuracy less). 
For MNIST with $\aleph =100$, the accuracy of instances drops slightly as the $k$ increases from 2 to 3, while using input features selected with $k=3$ leads to slightly higher accuracy for FashionMNIST with $\aleph = 100$. 
One reason behind this difference could be that input pixels of MNIST and FashionMNIST are activated in different patterns. In MNIST, there are 
more active pixels in the center, and the peripheral pixels stay inactive over the training data. 
Given the distribution constraints with $k=2$ (see Fig.~\ref{fig:mnist_100_k2}), the selected pixels stay away from the peripheral area. However, when $k=3$ (see Fig.~\ref{fig:mnist_100_k3}), more pixels are forced to be chosen from the upper right-hand and bottom left-hand corners. In contrast, the active pixels are more evenly spread across the full input in FashionMNIST. Hence, as $k$ increases (see Fig.~\ref{fig:fashionmnist_100_k2} and Fig.~\ref{fig:fashionmnist_100_k3}), the active pixels remain well-covered by the evenly scattered selected pixels.

Next, Table~\ref{tab:compare_k_mnist_dnn} and Table~\ref{tab:compare_k_fashionmnist_dnn} also compare OSIF using different DNN architectures, \textit{i.e.}, $2\times 10$ and $2 \times 20$. The accuracy is $94\%$ (resp. $84\%$) for the former and  $96\%$ (resp. $86\%$) for the latter for MNIST (resp. FashionMNIST).
The results show that even using a simple DNN for feature selection using OSIF can produce feature sets that achieve good accuracy, when appropriately large classification models are trained on the selected features. For MNIST (see Table~\ref{tab:compare_k_mnist_dnn}), the former model has accuracy at most 2 points worse than the latter model when $\aleph =50$. When $\aleph =100$, both models achieve similar levels of performance. As for Fashion, the difference is at most 1 point for $\aleph\in \{50,100\}$. Hence, we cannot observe a clear difference between the two OSIF DNN models in terms of feature selection quality.

Finally, we would like to make a brief remark on the improvement in accuracy by increasing the size of the DNNs for classification, \textit{i.e.,} from $2\times 20$  to $2\times 40$, for both MNIST and FashionMNIST. Unsurprisingly, using larger DNNs results in overall higher accuracy. More importantly, the performance of the proposed input feature selection seems to be robust toward the final architecture. For both architectures, we observe a similarly slight reduction in classification accuracy related to the reduction in number of input features (pixels).

\subsubsection{Comparisons between feature selection approaches}

In this section, we compare the performance of  the MILP-based features-selection approach (\textit{i.e.,} Alg.~\ref{alg:input_select}) to some other simple feature-selection approaches. 
The other feature selection techniques considered in the comparison are  random feature selection, data-based feature selection, and DNN weights-based. 
In the following paragraphs, we briefly describe the feature selection algorithms that are used as reference points for the comparison.

The random feature selection uniformly samples a subset of $\aleph $ input features, and we present the average accuracy with the standard deviation over 5 DNNs trained with input features independently selected by this approach. This approach is included, as it is the simplest approach to down-sample the full input.

The data-based feature selection is conducted in the following way: i) we first calculate the mean value of each input feature over the whole train dataset; ii) the features with the largest $\aleph$ mean are selected. The motivation behind this simple heuristic is that it selects the pixels that are most strongly colored over the training data, \textit{i.e.}, the strongest signals. For example, selecting a pixel uncolored in all images of the training data does not make sense as that input does not contain any information for the training data.

In the DNN weights-based approach, we use the same DNN models as we use in the MILP-based selection, but the inputs are now selected based on the weights of the inputs. For each input, we sum up the absolute values of all the weights from the input to the nodes in the consecutive layer and select the ones with the largest sum. This can be seen as a form of pruning of inputs, and the motivation is that inputs with small, or almost zero, weights should be less important as these inputs have less impact in the DNN.     

In Table~\ref{tab:compare_mnist_all}, we compare MILP-based feature selection (\textit{i.e.,} Alg.~\ref{alg:input_select}) with random selection, data-based selection, and weights-based selection. The figure for Alg.~\ref{alg:input_select} has the best result from Table~\ref{tab:compare_k_mnist_dnn}, where the accuracy of DNNs with sparse input features is only 5 points with $\aleph= 100$  less than the accuracy with full input features. It can be observed that our method has the best overall performance. For MNIST, the random selection has the worst performance, but the data-based selection and weights-based selection achieve a slightly worse performance than our method.

Table~\ref{tab:compare_fashionmnist_all} compares the performance of features selections on FashionMNIST. The results show a different pattern to MNIST. Our methods still have the best overall performance over different settings by maintaining the accuracy to 84\% (resp. 88\%) with $\aleph =100$ for DNNs $2\times 20$ (resp. $2\times 40$), while the accuracy of DNNs with full input features are $86\%$ (resp. $88\%$). While delivering the worst performance on MNIST, random selection has a very close performance to our method on FashionMNIST. The weights-based selection still lies in third overall, while the data-based selection is much worse than the other three methods (\textit{e.g.,} $58\%$ for the DNN $2\times 20$ with $\aleph=100$ and $62\%$ for the DNN $2\times 40$ with $\aleph =100$).
\begin{table}[!htbp]
    \centering
     \begin{tabular}{c|rrrrr|r}
      \hline
         \multirow{2}{*}{DNN} & \multicolumn{5}{c|}{Feature Selection Approaches} & \multirow{2}{*}{$\aleph = 784$} \\ \cline{2-6}
                            & $\aleph $ & MILP-based & Random & Data-based  &  Weights-based & \\ 
         \hline
         \multirow{2}{*}{$2\times 20$}   & 50    &  86.6\%   &   $77.4 \pm 2.2\%$  &   81.3\%  &   80.7\%    &     \multirow{2}{*}{95.7\%} \\
                                        & 100    &  91.2\%   &   $86.0\pm 2.5 \%$  &  89.2\%   &   89.0\%    &       \\
        \hline
         \multirow{2}{*}{$2\times 40$}  & 50    &  89.2\%    &   $76.0\pm 3.9\%$  &   85.1\%   & 83.3\%   & \multirow{2}{*}{97.1\%} \\
                                        & 100    &  93.4\%    &   $90.6 \pm 1.1\%$  &   92.4\%  &  91.2\%  & \\
        \hline
    \end{tabular}
    \caption{Accuracy of DNNs with sparse input features selected by different methods on MNIST}
    \label{tab:compare_mnist_all}
\end{table}
\begin{table}[!htbp]
    \centering
      \begin{tabular}{c|rrrrr|r}
      \hline
         \multirow{2}{*}{DNN} & \multicolumn{5}{c|}{Feature Selection Approaches} & \multirow{2}{*}{$\aleph = 784$} \\ \cline{2-6}
                            & $\aleph $ & MILP-based & Random & Data-based  &  Weights-based \\ 
         \hline
         \multirow{2}{*}{$2\times 20$}   & 50    &   77.9\%      &    $77.2\pm 0.8\%$       &   49.6\%   & 73.8\%  & \multirow{2}{*}{86.3\%}  \\
          & 100    &       82.3\%    &      $80.6\pm 1.1 \%$   &     58.4\%    & 80.3\%        \\
        \hline
         \multirow{2}{*}{$2\times 40$}     & 50   &    79.6\%     &             $78.8\pm 0.4 \%$    &    51.6\%   & 74.6\%  & \multirow{2}{*}{87.5\%} \\
          & 100    &    84.0\%       &     $82.8 \pm 0.4\%$          &       62.3\%     &   81.7\% &\\
        \hline
    \end{tabular}
    \caption{Accuracy of DNNs with sparse input features selected by different methods on FashionMNIST}
    \label{tab:compare_fashionmnist_all}
\end{table}

The weights-based selection performs decently on both data sets compared to random and data-based selection. However, based on the results it is clear that MILP-based selection (\textit{i.e.,} Alg.~\ref{alg:input_select}) can extract more knowledge from the DNN model regarding the importance of inputs compared to simply analyzing the weights. The overall performance of MILP-based selection is more stable than other feature selection methods on both datasets.

\subsection{Robustness to adversarial inputs}

The robustness of a trained DNN classifier can also be analyzed using MILP, \textit{e.g.,} in verification or finding adversarial input. We use the minimal distorted adversary as a measure of model robustness $x$ \lion{under $l_{\infty}$ norm~\cite{goodfellow2014explaining}.} For a given image $x_\mathrm{image}$, the minimal adversary problem \cite{tsay2021partition} can be formulated as~
\begin{equation}\label{eq:adversary}
    \begin{aligned}
        \min~ &\epsilon \\
      \text{s.t.} &\textrm{\eqref{eq:relu_1}--\eqref{eq:relu_var_size}},
         x_i^L \leq x_j^L, 
     || {x}^1 - {x}_\mathrm{image}||_\infty \leq \epsilon,
    \end{aligned}
\end{equation}
where $i$ is the true label of image $x_\mathrm{image}$ and $j$ is an adversarial label. 
Simply put, problem \eqref{eq:adversary} finds the smallest perturbation, defined by the $\ell_\infty$ norm, such that the trained DNN erroneously classifies image $x_\mathrm{image}$ as the adversarial label $j$.
We hypothesize that DNNs trained with fewer (well-selected) features are more robust to such attacks, as there are fewer inputs as degrees of freedom. 
Furthermore, we note that the robustness of smaller DNNs can be analyzed with significantly less computational effort. 
\begin{table}[!htbp]
    \centering
      \begin{tabular}{c|rrr|rrr}
      \hline
         \multirow{2}{*}{DNN}& \multicolumn{3}{c|}{MNIST} & \multicolumn{3}{c}{FashionMNIST}\\ 
         \cline{2-7}
          & $\aleph $ & $\epsilon$ ($\times 10^{-2}$) & Average $\Delta$ &  $\aleph $ & $\epsilon$ ($\times 10^{-2}$) & Average $\Delta$ \\ 
         \hline
         \multirow{3}{*}{$2\times 20$} & 50 & 16.1$\pm$7.9  & 69.9\% & 50 & 12.1$\pm$7.2 & 65.4\% \\
          & 100 & 14.0$\pm$6.1 & 33.6\% & 100 & 11.7$\pm$5.8 & 44.3\%     \\
          & 784 & 10.1$\pm$4.7 & - & 784 & 8.9$\pm$3.7 & - \\
        \hline
 
         \multirow{3}{*}{$2\times 40$} & 50 & 15.2$\pm$7.2 & 57.2\% & 50 & 12.9$\pm$7.6 & 89.1\% \\
          & 100    & 12.4$\pm$4.6 & 42.6\% & 100 & 11.5$\pm$5.5 & 52.7\% \\
          & 784    & 10.4$\pm$5.1 & - & 784 & 8.4$\pm$3.4 & - \\
        \hline
    \end{tabular}
    \caption{Minimal adversarial distance $\epsilon$ for trained DNNs.}
    \label{tab:compare_robustness}
\end{table}

Table~\ref{tab:compare_robustness} shows the minimal adversarial distance (mean and standard deviation over 100 instances), defined by \eqref{eq:adversary} for DNNs trained on MNIST and FashionMNIST with MILP-based feature selection. 
The adversaries are generated for the first 100 instances of the respective test datasets, with adversarial labels selected randomly. 
Furthermore, we report the mean percentage increase, $\Delta$, in minimal adversarial distance over the 100 instances for the reduced input DNNs compared to the full input. 
In all cases, reducing the the number of inputs $\aleph$ results in a more robust classifier. For the $2 \times 40$ DNN trained on FashionMNIST, reducing the number of inputs from 784 to 50 increases the mean minimal adversarial distance by almost 90\%, with a loss in accuracy of $<$10\%. 


\section{Conclusion} \label{sec:summary}
In the paper, we have presented an MILP-based framework using trained DNNs to extract information about salient input features. The proposed algorithm is able to drastically reduce the size of the input by using the input features that are most important for each category according to the DNN, given a regularization on the input size and spread of selected features. The numerical results show that the proposed algorithm is able to efficiently select a small set of features for which a good prediction accuracy can be obtained. The results also show that the proposed input feature selection can improve the robustness toward adversarial attacks.  
%
%
\bibliographystyle{splncs04}
\bibliography{Ref}
\end{document}